\documentclass[12pt]{amsart}

\usepackage{amsmath,amsthm,amsfonts,latexsym,amscd,amssymb,enumerate,times}
\usepackage{amscd}
\newtheorem{theorem}{Theorem}[section]
\newtheorem{thm}[theorem]{Theorem}
\newtheorem{prop}[theorem]{Proposition}
\newtheorem{lem}[theorem]{Lemma}
\newtheorem{cor}[theorem]{Corollary}
\newtheorem{conj}[theorem]{Conjecture}

\theoremstyle{remark}

\theoremstyle{definition}

\theoremstyle{definition}

\theoremstyle{remark}

\DeclareMathOperator{\im}{im}      

\DeclareMathOperator{\ind}{ind}

\newcommand{\rationals}{\mathbb{Q}}

\newcommand{\tensor}{\otimes}

\DeclareMathOperator{\End}{End}
\DeclareMathOperator{\ch}{{\rm ch}}

\newcommand{\K}{\mathbb{ K}}
\newcommand{\C}{\mathbb{ C}}

\newcommand{\N}{\mathbb{ N}}

\newcommand{\Z}{\mathbb{ Z}}

\newcommand{\Q}{\mathbb{ Q}}
\newcommand{\R}{\mathbb{ R}}

{\catcode`@=11\global\let\c@equation=\c@theorem}

\begin{document}

\title{The strong Novikov conjecture for low degree cohomology}

\author{Bernhard Hanke and Thomas Schick}

\begin{abstract} 
We show that for each discrete group $\Gamma$, the 
rational assembly map 
\[
    K_*(B\Gamma) \otimes \Q  \to K_*(C^*_{max} \Gamma) \otimes \Q
\]
is injective on classes dual to $\Lambda^* \subset H^*(B\Gamma;\Q)$, where 
$\Lambda^*$ is the subring generated by cohomology classes of degree 
at most $2$ (and where the pairing uses the Chern character). Our result
implies homotopy invariance of higher
signatures associated to classes in $\Lambda^*$. This consequence was   
first established by Connes-Gromov-Moscovici \cite{CGM} and Mathai \cite{Mat}. 

Our approach is based on the construction 
of flat twisting bundles out of sequences of almost flat 
bundles as first described in our work \cite{HS}. In contrast to 
the argument in \cite{Mat}, our approach is independent of 
(and indeed gives a new proof of) the result of Hilsum-Skandalis \cite{HSka} on 
the homotopy invariance of the index of the signature operator 
twisted with bundles of small curvature. 
\end{abstract} 

\address{Universit{\"a}t M{\"u}nchen \\ Germany \\}
\email{hanke@mathematik.uni-muenchen.de \vspace{1cm}}

\address{Universit{\"a}t G{\"o}ttingen  \\
Germany \\}
\email{schick@uni-math.gwdg.de}

\thanks{Both authors are members of the 
DFG emphasis programme ``Globale Differentialgeometrie'' whose 
support is gratefully acknowledged.}

\maketitle

\section{Introduction}

Throughout this paper, we use complex $K$-theory. Let $\Gamma$ be a discrete 
group and denote by $C^*_{max} \Gamma$ the maximal group $C^*$-algebra of $\Gamma$. 
Recall the following form of the strong Novikov conjecture.  

\begin{conj} The Baum-Connes assembly map 
\[
    A \colon  K_*(B\Gamma) \to K_*(C^*_{max} \Gamma) 
\]
is injective after tensoring with the rational numbers. 
\end{conj} 

The Chern character 
\[
    \ch \colon  K(-) \to H(-;\Q)
\]
is a natural transformation (of $\Z/2$-graded multiplicative (co-)homology 
theories) from 
$K$-homology to rational singular homology   (both theories being defined 
in the homotopy theoretic sense). Let 
\[
    \Lambda^*(\Gamma) \subset H^*(B\Gamma;\Q) 
\]
be the subring generated by classes of degree at most $2$. 

In this paper we verify the strong Novikov conjecture 
for classes dual to elements in $\Lambda^*(\Gamma)$. 

\begin{thm} \label{main} Let $h \in K_*(B\Gamma)$ be a $K$-homology class
such that the map  
\[
      \Lambda^*(\Gamma) \to \Q \, , ~ \gamma \mapsto \langle \gamma, \ch(h) \rangle \, , 
\]
given by the Kronecker pairing
is nonzero. Then 
\[
    A(h) \neq 0 \, . 
\]
\end{thm} 

As a corollary we obtain the following result on homotopy invariance of
higher signatures ($\mathcal{L}(M)$ denotes the $L$-polynomial of $M$).

\begin{cor}[\cite{CGM}, \cite{Mat}] Let $M$ be a closed connected oriented smooth manifold, let $\Gamma$ be a
discrete group and let $f \colon M  \to B\Gamma$ be a continuous map. Then 
for all $c \in \Lambda^*(\Gamma)$, the higher signature
\[
    \langle \mathcal{L}(M) \cup f^*(c) , [M] \rangle 
\]
is an oriented homotopy invariant. 
\end{cor} 

The discussion of this result in \cite{Mat} is based on a theorem of Hilsum-Skandalis \cite{HSka} saying 
that the index of the signature operator twisted with Hilbert $A$-module bundles 
($A$ being a $C^*$-algebra) of  
small curvature is an oriented homotopy invariant. For flat twisting bundles this result was 
known before (see e.g. \cite{Kam, Kas, Mis}). 
Our proof of Theorem \ref{main} is independent of \cite{HSka}. 
Indeed, we will illustrate in the last section of this paper how 
our methods allow the reduction of the  Hilsum-Skandalis theorem to the case of flat twisting bundles. 

\section{Proof of the main theorem} 

By a standard suspension argument we may assume that $h \in K_0(B\Gamma)$. Because each 
discrete group is the direct limit of finitely presented groups and because 
the classifying space construction, the formation of $C^*_{max}$ and the $K$-theory 
functors commute with direct limits, it is enough to treat the case of finitely 
presented $\Gamma$.

Due to the geometric description of $K$-homology by Baum-Douglas \cite{BD},
elaborated in \cite{BHS},  
there is triple $(M,E,\phi)$, where $M$ is 
an even dimensional closed connected spin-manifold, $E \to M$ is a virtual
(i.e.~$\Z/2$-graded) hermitian
complex vector bundle of finite dimension
and $\phi: M \to B\Gamma$ is a continuous map, such that 
\[
     \phi_*(E \cap [M]_K) = h \, . 
\]
(Strictly speaking, \cite{BHS} only provides a spin${}^c$-manifold; it is an
exercise to obtain a spin manifold such that the associated spin${}^c$-structure
represents the right $K$-homology class).
In the formula, we consider $E$ as a class in $K^0(M)$, and $[M]_K$ 
is the $K$-theoretic orientation class of the spin manifold $M$. 

As $\Gamma$ is finitely presented, we can and will assume
that the map $\phi: M \to B \pi_1(M)$ induces an isomorphism of fundamental 
groups.

Now chose $h \in K_0(B\Gamma)$ as in the main theorem and let 
$(M,E,\phi)$ be a triple representing $h$.
Let 
\[
    \nu = E \Gamma  \times_{\Gamma} C^*_{max} \Gamma 
\]
be the canonical flat $C^*_{max} \Gamma$-module bundle over $B\Gamma$. 

Denoting by $S^{\pm} \to M$ the bundles of positive or negative
complex spinors on $M$, respectively, 
the assembly map is described as follows. 

\begin{prop} The element 
\[
    A(h) \in K_0(C^*_{max} \Gamma)
\]
is equal to the Mishchenko-Fomenko index of 
\[
    D_{E \otimes \phi^*(\nu)} \colon \Gamma((S \otimes E)^+ \otimes \phi^*(\nu)) \to \Gamma((S \otimes E)^- \otimes \phi^*(\nu) ) \, , 
\]
the Dirac operator on $M$ twisted with the virtual bundle 
$E \otimes \phi^*(\nu)$. Here,  $E$ is equipped with an arbitrary hermitian connection.  
\end{prop} 

This description of the assembly map will be used in order to show that 
$A(h) \neq 0$. 

For that purpose, let $c \in \Lambda^*(\Gamma)\subset H^*(B\Gamma,\rationals)$ be such that 
\[
      \langle c , \ch(h)  \rangle \neq 0 \in \Q \, . 
\] 
We can assume that $c \in H^*(B\Gamma; \Z)$. In order to keep the exposition transparent, 
let us first assume that $c \in H^2(B\Gamma;\Z)$.

Let $L \to B\Gamma$ be a complex hermitian line bundle classified by $c$. We
pick  
a unitary connection on the pull back bundle $L'=\phi^*(L) \to M$ and denote by 
$\omega \in \Omega^2(M;i\R)$ its curvature form. 
Let $\pi : \widetilde M \to M$ be the universal cover. Because the universal cover of $B\Gamma$ is contractible,
the bundle $\pi^*(L') \to \widetilde M$ is trivial. Fix a unitary trivialization  
$\pi^*(L') \cong \widetilde M \times \C$. 
With respect to this trivialization, the induced connection on $\pi^*(L')$ 
is given by a $1$-form $\eta \in \Omega^1(\widetilde M; i \R)$. 

Because $U(1)$ is abelian, the curvature form of this connection is equal 
to $ d \eta$ which in turn coincides with $\pi^*(\omega)$ by naturality. However, contrary to 
the form $\pi^*(\omega)$, the connection form $\eta$ is in general not invariant under 
the deck transformation group (this would imply that  $\omega$ represents the zero 
class in $H^2(M)$). 

We will use the bundle $\pi^*(L')$ in order to construct a flat 
$A$-module bundle $W \to M$ with an appropriate $C^*$-algebra $A$
along the lines of \cite{HS}. 
The flat bundle $W$ will induce a holonomy representation 
\[
   C^*_{max} \pi_1(M)   \to A 
\]
whose induced map in $K$-theory will be used to detect non-triviality 
of the element $A(h)$ appearing in the main theorem.

The details are as follows. For  
$t \in [0,1]$ we consider the connection on $\widetilde M \times
\C$ associated to the $1$-form  $t \cdot \eta$. 
The corresponding curvature form is equal to $t \cdot \pi^*(\omega)$ and
this is invariant under deck transformations (note that this is in general 
not true for the forms $t \cdot \eta$ if $t \neq 0$).

We would like to use the bundle $L'$ to construct a family of almost flat bundles 
$(P_t)_{t \in [0,1]}$ (cf. \cite{HS}, Section 2)  so that 
\cite[Theorem 2.1]{HS} 
can be applied in order to obtain an ``infinite product bundle''
\[   
     V = \prod_{n \in \N} P_{1/n} \to M  
\]
whose quotient by the corresponding infinite sum bundle will be the desired flat bundle $W \to M$, 
cf. \cite[Proposition 3.4]{HS}. However, by Chern-Weil theory
it is in general impossible to produce a finite dimensional bundle on $M$ whose 
curvature form is equal to $t \cdot \omega$, $0 < t < 1$ (the associated 
Chern class would not be integral). 

We bypass this difficulty by allowing infinite dimensional 
bundles. Consider the Hilbert space bundle 
\[
             \mu =  \widetilde M \times_{\Gamma} l^2(\Gamma) \to M \, .  
\]
Here, $\Gamma$ acts on the left of $l^2(\Gamma)$ by the formula
\[
     (\gamma \psi)(x) := \psi(x \cdot \gamma^{-1}) \, , ~ x, \gamma \in \Gamma \, . 
\]
The forms $t \eta \in \Omega^1(\widetilde M)$ induce a family 
of connections $\nabla^t$ on $\mu$, the connection  $\nabla^t$ being induced 
by the $\Gamma$-invariant connection on $\widetilde M \times l^2(\Gamma)$   
which on the subbundle
\[
    \widetilde M \times \C \cdot 1_{\gamma} \subset \widetilde M \times l^2(\Gamma) 
\]
(identified canonically with $\widetilde M \times \C$) coincides with 
$(\gamma^{-1})^*(t \eta)$. 

We wish to regard the bundles $(\mu, \nabla^t)$ as twisting bundles for the Dirac 
operator on $M$. The index of this Dirac operator will live in the $K$-theory 
of an appropriate $C^*$-algebra $A_t$ that we think of as the holonomy algebra of $\mu$ 
with respect to the connection $\nabla^t$.

To define these algebras, let us choose a base point $p \in M$ 
and a point $q \in \widetilde M$ above $p$. 
We identify the 
fibre over $p$ with the Hilbert space $l^2(\Gamma)$. 

Now let
\[
      A_t \subset B(l^2(\Gamma)) 
\]
be the norm-linear closure of all maps $l^2(\Gamma) \to l^2(\Gamma)$ that
arise from parallel transport with respect to $\nabla^t$
along closed curves in $M$ based at $p$. Further, we define a bundle $P_t \to M$ whose 
fibre over $x \in M$ is given by the norm-linear closure (in ${\rm Hom}(\mu_p, \mu_x)$) of  
all isomorphisms  $\mu_p \to \mu_x$
arising from parallel transport with respect to $\nabla^t$  along  smooth curves connecting $p$ with $x$. 
In this way we obtain, for each $t \in [0,1]$, a locally trivial bundle $P_t$ consisting of free $A_t$-modules of rank 
$1$ (the $A_t$-module structure given 
by precomposition) and equipped with $A_t$-linear connections. For the notions relevant 
in this context, we refer the reader e.g.~to \cite{Sch}. Parallel transport on $P_t$ 
along a curve connecting $x$ with $x'$ is induced by parallel transport $\mu_x \to \mu_{x'}$ with 
respect to $\nabla^t$. The bundle $\mu$ may be recovered from the ''principal 
bundle`` $P_t$ as an associated bundle, i.e. $\mu = P_t \times_{A_t} l^2(\Gamma)$. 

The next lemma is crucial for the calculations which will follow.

\begin{lem} Each of the algebras $A_t$ carries a canonical trace 
\[
     \tau_t : A_t \to \C
\]
given by 
\[
      \tau_t(\psi) :=  \langle \psi(1_e), 1_e \rangle 
\]
where $1_e \in l^2(\Gamma)$ is the characteristic function of the neutral element and $\langle -, - \rangle$ 
is the inner product on $l^2(\Gamma)$. 
\end{lem} 

\begin{proof} Let $\gamma$ and $\gamma'$ be two closed smooth curves based 
at $p \in M$ and let $\phi_{\gamma}$ and $\phi_{\gamma'}$ be parallel 
transport along $\gamma$ and $\gamma'$. We show that 
\[
   \tau_t(\phi_{\gamma} \cdot \phi_{\gamma'}) = \tau_t(\phi_{\gamma'} \cdot \phi_{\gamma}) \, . 
\]
We will assume from now on that $\gamma$ and $\gamma'$ represent elements in $\pi_1(M,p)$ 
that are inverse to each other (otherwise, both sides of the above equation are zero). 
We lift the composed 
curves $\gamma \cdot \gamma'$ and $\gamma' \cdot \gamma$ to curves $\chi$ and $\chi'$ starting at 
$q$ in the universal cover $\widetilde M$. By assumption, both $\chi$ and $\chi'$ are closed curves. 
We need to compare parallel transport of the element 
\[
      (q,1) \in \widetilde M \times \C 
\]
along $\chi$ and $\chi'$. Denoting the result of parallel transport along these 
curves by $(q,\xi_{\chi})$ and $(q, \xi_{\chi'})$, respectively, we have 
\[
    \xi_{\chi} = \exp \big( \int_{[0,1]} - t\eta_{\chi(\tau)}(\dot{\chi}(\tau)) d\tau \big) \, , 
\]
with the exponential map on the Lie group $S^1$. By use of Stoke's formula, the last expression 
is equal to 
\[
    \exp \big( \int_{D} -d (t\eta) dx \big) = \exp \big( \int_D -t \pi^*(\omega)  dx \big) \, , 
\]
where $D : D^2 \to \widetilde M$ is a disk with boundary $\chi$ (recall that $\widetilde M$ is simply connected). 
However, the curve $\chi'$ is obtained (up to reparametrization) from $\chi$ by applying the 
deck transformation corresponding to $\gamma'$. This implies - grace to the invariance of 
the form $\pi^*(\omega)$ under deck transformation - that the last expression is equal 
to $\xi_{\chi'}$. Using continuity and linearity of $\tau_t$, this shows our
assertion.
\end{proof}

We now equip $E$ with a unitary connection and consider the twisted 
Dirac operator 
\[
    D_{E \otimes P_t} : \Gamma((S \otimes E)^+ \otimes P_t) \to \Gamma((S \otimes E)^- \otimes P_t) \, . 
\]
The index $\ind(D_{E \otimes P_t}) \in K_0(A_t)$ satisfies  (up to sign) the equation 
\[
 \tau_t(\ind(D_{E \otimes P_t})) = \langle \ch(\mathcal{A}(M)) \cup \ch(E) \cup \ch_{\tau_t}(P_t) , [M] \rangle 
\]
by the Mishchenko-Fomenko index theorem (see \cite[Theorem 6.9]{Sch}). 
Because the Poincar\'e dual of $\ch(\mathcal{A}(M))$ is equal to $\ch([M]_K)$, the last expression 
equals  $\langle  \ch_{\tau_t}(P_t) ,  \ch([E] \cap [M]_K) \rangle$. Recall
that the choice of $(M,E,\phi)$ 
implies  $\phi_*\ch([E]\cap [M]_K) = h$.
The construction of the connection  $\nabla^t$ and the definition of $\ch_{\tau_t}$ (see \cite[Definition 5.1.]{Sch}) show 
\[
     \ch_{\tau_t}(P_t) = \exp(t\phi^*c) \in H^*(M ; \R)   
\]
so that finally 
\[
   \tau_{t}(\ind(D_{E \otimes P_t})) = \langle \exp(t \phi^*c)  , \ch([E]\cap
   [M]_K) \rangle = \langle \exp(tc) , \ch(h)\rangle\in \R[t] \, ,  
\]
a polynomial which is different from zero by the assumption $\langle c, \ch(h) \rangle \neq 0$. 

We wish to use this calculation in order to detect non-triviality of $A(h) \in K_0(C^*_{max} \Gamma)$. 
This is done by constructing a flat bundle on $M$ out of the sequence of 
bundles $(P_{1/n})_{n \in \N}$ with connections. This sequence 
is almost flat in the sense of \cite[Section 2]{HS}. Therefore, applying \cite[Theorem 2.1.]{HS} 
we obtain a smooth ``infinite product bundle''
\[
    V = \prod_{n \in \N} P_{1/n} \to M 
\]
equipped with a connection that we may think of as the product of the connections on $P_{1/n}$. 
After passing to the quotient by the infinite sum of the bundles $P_{1/n}$, 
we obtain a smooth bundle 
\[
    W = \big( \prod P_{1/n} \big) / \big( \bigoplus P_{1/n} \big) \to M 
\]
of free modules of rank $1$ over the $C^*$-algebra 
\[
    A = \big( \prod A_{1/n} \big) / \big( \overline{\bigoplus A_{1/n}} \big)  
\]
which carries an induced flat connection. For more details of this discussion, we 
refer the reader to \cite{HS}, especially to Section 2 and the statements before Proposition 3.4. 

The flat bundle $W$ induces a unitary holonomy representation of $\pi_1(M)$. Because 
$\phi$ induces an isomorphism $\pi_1(M) \cong \Gamma$ by the choice of the triple $(M,\phi,E)$, 
we hence obtain a $C^*$-algebra map 
\[
   \psi \colon C^*_{max} \Gamma \to A  
\]
using the universal property of $C^*_{max}$. The induced map in $K$-theory can be analyzed 
in terms of the  $KK$-theoretic description of the index map (cf. \cite[Lemma 3.1]{HS}). One 
concludes that $\psi_*(A(h)) \in K_0(A)$ is equal to the index of the twisted Dirac  operator 
\[
     D_{E \otimes W} \colon \Gamma((S \otimes E)^+ \otimes W) \to \Gamma((S \otimes E)^- \otimes W)   
\]
and this in turn equals the image of the index of the twisted Dirac operator 
\[  
     D_{E \otimes V} \colon \Gamma((S \otimes E)^+ \otimes V) \to \Gamma((S \otimes E)^- \otimes V) 
\]
under the canonical map $p_*  \colon K_0(\prod A_{1/n}) \to K_0(A)$ induced by the canonical 
projection $p \colon \prod A_{1/n} \to A$. 

We have a short exact sequence 
\[
    K_0( \bigoplus  A_{1/n}) \to K_0(\prod  A_{1/n}) \stackrel{p_*}{\to} K_0(A)   
\]
where the left hand group is canonically isomorphic to the algebraic direct sum 
$\bigoplus_{n\in \N} K_0(A_{1/n})$. Furthermore, the traces $\tau_{1/n} \colon A_{1/n} \to \C$ 
all have norm $1$ and hence induce a trace 
\[    
      \prod A_{1/n} \to \prod \C \, . 
\]
Using the canonical isomorphism $K_0( \bigoplus A_{1/n} ) \cong \bigoplus K_0(A_{1/n})$, 
we finally get a trace map 
\[ 
   \tau \colon \im p_* \to \big( \prod \C \big) / \big( \bigoplus \C  \big) \, . 
\]
Note that the direct sum on the right is understood in the algebraic sense: 
Any element has only finitely many components different from $0$.  Assuming 
$\ind(D_{E \otimes W}) = 0$ we conclude that the element in $\prod\C / \bigoplus \C$ 
represented by 
\[
     (\ind(D_{E \otimes P_{1/n}}))_{n \in \N} 
\]
is equal to zero. Combining this with our previous calculation, this means that the polynomial 
\[
     \langle \exp(tc) , \ch(h) \rangle \in \R[t] 
\]
is equal to zero for all but finitely many values $t = 1/n$, $n \in \N$. This implies 
that this  polynomial is identically  zero and hence in particular 
\[
    \langle c, \ch(h) \rangle  = 0 
\]
in contradiction to our assumption. The proof of the main 
theorem is therefore complete, if $c \in H^2(B\Gamma;\Z)$. 

We will now discuss the case of general $c \in \Lambda^*(\Gamma)$ (still assuming $h \in K_0(B\Gamma)$). 
If  
\[
    c = c_1 \cup \ldots \cup c_k 
\]
is a product of classes in $H^2(B\Gamma;\Z)$,  we replace the bundle $L$ in the above 
argument by the tensor product bundle 
\[
    L :=  L_1 \otimes L_2 \otimes \ldots \otimes L_k \to B \Gamma 
\]
where the line bundle $L_i \to B\Gamma$ is classified by $c_i$. In an 
analogous fashion as before, we get bundles 
\[
   P_{t_1, \ldots, t_k} \to M 
\]
of Hilbert-$A_{t_1, t_2, \ldots, t_k}$-modules working with the connection 
\[
      t_1 \eta(1) + t_2 \eta(2) + \ldots + t_k \eta(k) 
\]
on the bundle $\widetilde M \times \C$ where $\eta(i)$ is a connection induced
from $L_i$ and where each $t_i \in [0,1]$. 

Assuming that $A(h) = 0$ we can conclude in this case that the 
polynomial 
\[
    \langle  \exp(t_1  c_1) \cdot \ldots \cdot \exp (t_k c_k) , \ch(h) \rangle \in \R[t_1, \ldots, t_k] 
\]
is equal to zero for all but finitely many 
\[
    (t_1, \ldots, t_k ) = (1/n_1, \ldots, 1/n_k) 
\]
where $n_i \in \N$. Hence this polynomial is identically $0$ and in particular 
\[
         \langle c_1 \cup \ldots \cup c_k, \ch(h) \rangle  = 0 
\]
again contradicting our assumption. In the most general case, $c$ being 
a sum 
\[
    c = c(1) + \ldots + c(k) 
\]
where each $c(i)$ is a product of two dimensional classes in $H^2(B\Gamma;\Z)$, the 
assumption 
\[
   \langle c , \ch(h) \rangle \neq 0 
\]
implies that already for one summand $c(i)$, we have $\langle c(i) , \ch(h) \rangle \neq 0$ 
so that we are reduced to the previous case. 

Because we already used a suspension argument in order to restrict attention to 
classes in $K_0(B\Gamma)$, this finishes the proof of the main theorem. 

\section{Hilsum-Skandalis revisited} 

If $M$ is an oriented Riemannian manifold of even dimension, recall that $d+d^*$, the sum of the exterior 
differential and its formal adjoint, defines the signature operator 
\[
    D^{sign} \colon \Gamma(\Lambda_+^*(M)) \to \Gamma(\Lambda_-^*(M) ) \, , 
\]
where the $\pm$-signs indicate $\pm 1$-eigenspaces of the Hodge star operator. 

If $E \to M$ is a Hilbert $A$-module bundle over $M$, where $A$ is some $C^*$-algebra, we 
obtain the twisted signature operator 
\[
    D^{sign}_{E} \colon \Gamma(\Lambda_+^*(M)\tensor E) \to
    \Gamma(\Lambda_-^*(M)\tensor E) 
\]
which has an index $\ind{D^{sign}_E} \in K_0(A)$. The following theorem says that this class is 
an oriented homotopy invariant, if $E$ has small curvature.

\begin{thm}[Hilsum-Skandalis \cite{HSka}] \label{HilsumSka} Let $M$ and $M'$
  be closed oriented Riemannian manifolds of the same dimension  
and let $h \colon  M' \to M$ be an orientation preserving homotopy
equivalence. Then there exists a constant $c > 0$ 
with the following property: If $E \to M$ is a Hilbert $A$-module bundle
with connection $\nabla$ so that the associated curvature form
$\Omega_{\nabla} \in \Omega^2(M; \End{E})$ satisfies
the bound 
\[
      \| \Omega_{\nabla} \| < c \, ,  
\]
(the norm being defined by  the maximum norm on the unit sphere bundle in $\Lambda^2 M$ and 
the operator norm on each fibre $\End(E_x, E_x)$), then we have
\[
     \ind(D^{sign}_{f^*(E)}) = \ind(D^{sign}_{E}) \, . 
\]
\end{thm}

If $E$ is flat, i.e. $\Omega_{\nabla} = 0$, this result was proved in \cite{Kam,Kas, Mis}. 
In this section, we will explain briefly how this special case implies the general statement of Theorem \ref{HilsumSka}. 
Our argument is again based on the construction of a flat bundle out of a sequence of 
almost flat bundles \cite{HS}. 
 
Assuming that no $c$ with the stated property exists, we find 
a sequence of $C^*$-algebras $A_n$ and a sequence of Hilbert $A_n$-module bundles $E_n \to M$ 
with connections $\nabla_n$ so that $\|\Omega_{\nabla_n} \| < \frac{1}{n}$, but 
\[
    \ind(D^{sign}_{f^*(E_n)}) \neq  \ind(D^{sign}_{E_n})
\]
for all $n$. We obtain an almost flat sequence 
\[
    f^*(E_n) \cup  E_n \to M' \cup M 
\]
of Hilbert $A_n$-module bundles in the sense of \cite[Section 2]{HS} over the disjoint
 union $M' \cup M$. 
Applying again the methods in \cite{HS}, we obtain a flat Hilbert 
$A$-module bundle  
\[
     W = \big( \prod f^*(E_n) \cup \prod E_n \big)   / \big( \bigoplus f^*( E_n)  \cup \bigoplus E_n \big)   \to M' \cup M 
\]
where 
\[
    A = \big( \prod A_n \big) / \big( \overline{\bigoplus A_n}  \big) \, . 
\] 
The index $\ind(D^{sign}_W)$ of the signature operator on $(-M') \cup M$ (the minus-sign 
indicates reversal of orientation) twisted by $W$, vanishes by the results in \cite{Kam, Kas,Mis}. 

This conclusion leads to the following contradiction. 
We have a canonical isomorphism 
\[
    K_0(\bigoplus A_n) \cong \bigoplus K_0(A_n ) 
\]
and - assuming that each $A_n$ is unital and stable, the last property easily being 
achieved by replacing each $A_n$ by $A_n \otimes \K(l^2(\N))$ - a canonical 
isomorphism 
\[
    K_0(\prod A_n)  \cong  \prod K_0(A_n) \, , 
\]
compare the proof of \cite[Proposition 3.6]{HS}. The signature operator 
on $- M' \cup M$ twisted with the non-flat Hilbert $\prod A_n$-module bundle
\[
     V := \prod f^*(E_n) \cup \prod E_n 
\]
has an index 
\[
    \ind(D^{sign}_V) \in K_0(\prod A_n) \cong \prod K_0(A_n)
\]
which is different from $0$ for infinitely many factors by our assumption on the bundles $E_n$. 
On the other hand, under the canonical map 
\[
     p_* \colon  K_0(\prod A_n) \to K_0(A) 
\]
this index maps to the index of the signature operator twisted with the bundle $W$ and 
this index was identified  as $0$. Now the contradiction arises from 
the fact that by the long exact $K$-theory sequence, the kernel of $p_*$ 
is equal to $K_0(\bigoplus A_n) = \bigoplus K_0(A_n)$ (algebraic direct sum) 
and therefore $\ind(D^{sign}_V)$ 
is not contained in it.


\begin{thebibliography}{AAAA}

\bibitem{B} B. Blackadar, {\it $K$-theory for operator algebras},
Second edition, Cambridge University Press, Cambridge, 1998. 


\bibitem{BD} Baum, Paul and Douglas, Ronald G, {\it $K$-homology and index
    theory}.  Operator algebras and applications, Part I (Kingston, Ont.,
  1980),  117--173, Proc. Sympos. Pure Math., 38, Amer. Math. Soc.,
  Providence, R.I. (1982) 

\bibitem{BHS} Paul Baum, Nigel Higson, and Thomas Schick, {\it {{On the
        Equivalence of Geometric and Analytic K-Homology}}},  {Pure and Applied Mathematics Quarterly 3, 1--24, (2007)},
 {arXiv:math.KT/0701484}

\bibitem{CGM} A. Connes, M. Gromov, H. Moscovici, {\em Group cohomology with Lipschitz control and 
higher signatures}, Geom. Funct. Anal. {\bf 3}(1) (1993), 1-78



\bibitem{HS} B. Hanke, T. Schick, {\em Enlargeability and index theory}, J. Differential Geom. 
{\bf 74} (2006), 293-320. 

\bibitem{HSka} M. Hilsum, G. Skandalis, {\em Invariance par homotopie de la signature \`a 
coefficients dans un fibr\'e presque plat}, J. Reine Angew. Math {\bf 423} (1992), 73-99. 

\bibitem{Kam} J. Kaminker, J.G. Miller, {\em Homotopy invariance of the analytic signature operators
over $C^*$-algebras}, J. Op. Th. {\bf 14} (1985), 113-127.


\bibitem{Kas} G. Kasparov, {\em Topological invariance of elliptic operators I: $K$-homology}, Izv. Akad. Nauk. 
S.S.S.R. Ser. Mat. {\bf 39} (1975), 796-838; English transl. in Math. U.S.S.R. Izv. {\bf 9} (1975), 751-792. 

\bibitem{Mat} V. Mathai, {\em The Novikov conjecture for low degree cohomology classes}, 
Geom. Dedicata {\bf 99} (2003), 1-15.  

\bibitem{Mis} A. Mishchenko, {\em $C^*$-algebras and $K$-theory}, Lect. Notes Math. {\bf 763} (1979), 262-274. 


\bibitem{Sch} T. Schick, {\em $L^2$-index theorems, $KK$-theory
and connections}, New York J. Math. {\bf 11} (2005), 387-443. 



\end{thebibliography}
\end{document}